\input ./style/arxiv-general.cfg
\documentclass[aop,MSNbibl,dvips]{arximspdf}
\makeatletter
   \@ifpackageloaded{graphicx}{}{\usepackage{graphicx}}
\makeatother
\usepackage{mathrsfs}

\doi{10.1214/14-AOP967} 
\volume{43}
\issue{6}
\pubyear{2015}
\firstpage{3494}
\lastpage{3513}
\docsubty{FLA}

\makeatletter
\newproclaim{remark}{Remark}
\mathchardef\varDelta="0101
\def\cal{\mathscr}

\newcommand{\spe}{{\cal S}}
\newcommand{\PP}{{\cal P}}
\newcommand{\LL}{{\cal L}}
\newcommand{\WW}{{\cal W}}
\newcommand{\Ind}{{\mathrm{I}}}
\newcommand{\eps}{{\varepsilon}}
\newcommand{\ch}{\mathrm{ch}}
\newcommand{\smsp}{}
\newcommand{\e}{\mathbb{E}}
\newcommand{\q}{\mathbb{Q}}
\newcommand{\Natural}{\mathbb{N}}
\newcommand{\la}{\langle}
\newcommand{\ra}{\rangle}

\newtheorem{lemma}{Lemma}
\newtheorem{theorem}{Theorem}
\newtheorem{corollary}{Corollary}

\makeatother

\begin{document}
\begin{frontmatter}

\title{The free energy in a multi-species Sherrington--Kirkpatrick model}
\runtitle{Multi-species SK model}

\begin{aug}
\author{\fnms{Dmitry} \snm{Panchenko}\corref{}\ead[label=e1]{panchenk@math.toronto.edu}}
\address{Department of Mathematics\\
University of Toronto\\
Bahen Centre, 40 St. George St.\\
Toronto, Ontario M5S 2E4\\
Canada\\
\printead{e1}}
\runauthor{D.~Panchenko}
\affiliation{University of Toronto}
\end{aug}

\received{\smonth{2} \syear{2014}}
\revised{\smonth{8} \syear{2014}}

%
\begin{abstract}
The authors of
[\textit{Ann. Henri Poincar\'{e}}  \textbf{16} (2015) 691--708]
introduced a~multi-species version of the
Sherrington--Kirkpatrick model and suggested the analogue of the Parisi
formula for the free energy. Using a variant of Guerra's replica
symmetry breaking interpolation, they showed that, under certain
assumption on the interactions, the formula gives an upper bound on the
limit of the free energy. In this paper we prove that the bound is
sharp. This is achieved by developing a new multi-species form of the
Ghirlanda--Guerra identities and showing that they force the overlaps
within species to be completely determined by the overlaps of the whole system.
\end{abstract}

%

%
\end{frontmatter}
%
\section{Introduction and main results}
Recently, the following modification of the Sherrington--Kirkpatrick
model \cite{SK} was introduced in \cite{MS}. Given $N\geq1$, let us
denote by
%
\begin{equation}
\sigma= (\sigma_1,\ldots,\sigma_N) \in
\Sigma_N = \{-1,+1 \}^N
\end{equation}
a configuration of $N$ Ising spins. Consider a finite set $\spe$ that
will be fixed throughout the paper and, in particular, it does not
change with $N$. We emphasize this because we will often omit the
dependence of other objects on $N$. The elements of $\spe$ will be
called species and will be denoted by $s$ or $t$. Let us divide all
spin indices into disjoint groups indexed by the species
%
\begin{equation}
I = \{1,\ldots, N \} = \bigcup_{s\in\spe}
I_s. \label{species}
\end{equation}
These sets will, obviously, vary with $N$, and we will assume that
their cardinalities $N_s =|I_s|$ satisfy
%
\begin{equation}
\lim_{N\to\infty} \frac{N_s}{N} = \lambda_s \in(0,1)
\qquad\mbox{for all } s\in\spe.
\end{equation}
For simplicity of notation, we will omit the dependence of $\lambda_s^N
:= N_s/N$ on $N$ and will simply write $\lambda_s$. The Hamiltonian
proposed in \cite{MS} resembles the usual SK Hamiltonian,
%
\begin{equation}
H_N(\sigma) = \frac{1}{\sqrt{N}} \sum_{i,j =1}^N
g_{ij}\sigma_i \sigma_j, \label{SKH}
\end{equation}
where the interaction parameters $(g_{ij})$ are independent Gaussian
random variables, only now they are not necessarily identically
distributed but, instead, satisfy
%
\begin{equation}
\e g_{ij}^2 = \Delta_{st}^2 \qquad\mbox{if } i\in I_s, j\in I_t
\mbox{ for } s,t\in\spe.
\label{Delta}
\end{equation}
In other words, the variance of the interaction between $i$ and $j$
depends only on the species they belong to. We will make the same
assumptions on the matrix $\Delta^2 = (\Delta_{st}^2)_{s,t\in\spe}$ as
in \cite{MS}, namely, that it is symmetric and nonnegative definite,
%
\begin{equation}
\Delta_{st}^2 = \Delta_{ts}^2 \qquad\mbox{for all } s,t\in\spe\mbox{ and } \Delta^2 \geq0. \label{assumption}
\end{equation}
Let us denote the overlap of the restrictions of two spin
configurations to a given species $s\in\spe$ by
%
\begin{equation}
R_s\bigl(\sigma^1,\sigma^2\bigr) =
\frac{1}{N_s} \sum_{i\in I_s} \sigma_i^1
\sigma_i^2. \label{overlap}
\end{equation}
Then it is easy to see that the covariance of the Gaussian Hamiltonian
(\ref{SKH}) is given by
%
\begin{equation}
\frac{1}{N} \e H_N\bigl(\sigma^1
\bigr)H_N\bigl(\sigma^2\bigr) = \sum
_{s,t\in\spe} \Delta_{st}^2
\lambda_s \lambda_t R_s\bigl(\sigma
^1,\sigma ^2\bigr) R_t\bigl(
\sigma^1,\sigma^2\bigr). \label{Cov}
\end{equation}
This already gives some idea about the main new difficulty one
encounters in this model compared to the classical
Sherrington--Kirkpatrick model. Namely, now we will need to understand
the joint distributions of the overlap arrays in the thermodynamic
limit simultaneously for all species $s\in\spe$. Our main goal will be
to compute the limit of the free energy in this model,
%
\begin{equation}
F_N = \frac{1}{N}\smsp\e\log Z_N,\qquad \mbox{where }
Z_N = \sum_{\sigma\in\Sigma_N } \exp H_N(
\sigma). \label{FE}
\end{equation}
Notice that we do not consider the inverse temperature parameter here,
because it can be absorbed into the definition of the matrix $\Delta
^2$. One can also consider the externals fields that depend only on the
species but, since it does not affect any arguments in the paper, for
simplicity of notation we will omit them.

Under assumption (\ref{assumption}), the authors in \cite{MS} proved,
using the Guerra--Toninelli interpolation \cite{GuerraToninelli}, that
the free energy has a limit. They also proposed the following analogue
of the Parisi formula \cite{Parisi79,Parisi} for the free energy,
which was proved for the original SK model by Talagrand in \cite{TPF};
see also \cite{SG2}. Given integer $r\geq1$, consider a sequence
%
\begin{equation}
0=\zeta_{-1}< \zeta_0 <\cdots< \zeta_{r-1} <
\zeta_r = 1 \label{zetas}
\end{equation}
and, for each $s\in\spe$, a sequence
%
\begin{equation}
0=q^s_0\leq q^s_1 \leq\cdots
\leq q^s_{r-1}\leq q^s_r =1.
\label{qs}
\end{equation}
We will also consider two types of nondecreasing combinations of these
sequences as follows. For $0\leq\ell\leq r$, we define
%
\begin{equation}
Q_\ell= \sum_{s,t\in\spe} \Delta_{st}^2
\lambda_s \lambda_t q^s_\ell
q^t_\ell \quad\mbox{and}\quad Q^s_\ell= 2\sum
_{t\in\spe} \Delta_{st}^2
\lambda_t q^t_\ell \mbox{ for } s\in\spe.
\label{Qs}
\end{equation}
The meaning of these definitions will become clear when we look at the
covariance of the cavity fields in the Aizenman--Sims--Starr scheme in
Section~\ref{Sec5label}. Given these sequences, let us consider i.i.d.
standard Gaussian random variables $(\eta_\ell)_{1\leq\ell\leq r}$
and, for $s\in\spe$, define
%
\begin{equation}
X^s_r = \log\ch\smsp \sum_{1\leq\ell\leq r}
\eta_\ell \bigl(Q^s_{\ell} - Q^s_{\ell
-1}
\bigr)^{1/2}. \label{Xr}
\end{equation}
Recursively over $0\leq\ell\leq r-1$, we define
%
\begin{equation}
X^s_\ell=\frac{1}{\zeta_\ell}\log\e_\ell\exp
\zeta_\ell X^s_{\ell+1}, \label{Xelll}
\end{equation}
where $\e_\ell$ denotes the expectation with respect to $\eta_{\ell+1}$
only. Notice that $X^s_0$ are nonrandom. Finally, we define the
analogue of the Parisi functional by
%
\begin{equation}
\PP(\zeta,q) = \log2 + \sum_{s\in\spe}
\lambda_s X^s_0 -\frac{1}{2} \sum
_{0\leq\ell\leq r-1} \zeta_\ell (Q_{\ell+1} -
Q_{\ell
} ). \label{Pzeta}
\end{equation}
The main result of the paper is the following.

\begin{theorem}\label{ThFE}
Under the assumption (\ref{assumption}), the limit of the free energy
is given by
%
\begin{equation}
\lim_{N\to\infty} F_N = \inf\PP(\zeta,q), \label{Parisi}
\end{equation}
where the infimum is taken over $r\geq1$ and the sequences (\ref
{zetas}) and (\ref{qs}).
\end{theorem}

In \cite{MS}, the inequality $F_N \leq\inf\PP(\zeta,q)$ was proved
under assumption (\ref{assumption}) using the analogue of Guerra's
replica symmetry breaking interpolation \cite{Guerra}. For convenience,
we will reproduce this result in Section~\ref{Sec2label} in the
formalism of the Ruelle probability cascades, which will also allow us
to introduce several objects that will be used in the subsequent
sections. In this paper we will prove the matching lower bound using
the analogue of the Aizenman--Sims--Starr scheme \cite{AS2} and, in
this part, the assumption $\Delta^2\geq0$ will not be needed. The
approach was applied previously in various situations in \cite{PPF} and
\cite{WKC} and is based on the ultrametricity result in~\cite{PUltra}.
As we mentioned above, in the multi-species model we encounter a new
nontrivial obstacle. Namely, we need to describe the joint
distribution of the overlap arrays simultaneously for all species, and
even though it is clear that the marginal distribution of each array
will be generated by the Ruelle probability cascades as in the SK
model, it is not at all clear what their joint distribution should be.
We will develop an approach to overcome this obstacle in Sections~\ref{Sec3label} and \ref{Sec4label}. In Section~\ref{Sec3label} we will
prove a multi-species version of the Ghirlanda--Guerra identities,
which are similar to the original Ghirlanda--Guerra identities \cite{GG}, but apply to generic overlaps that may depend on the overlaps of
all species. Using these identities, we will show in Section~\ref{Sec4label} that the overlaps of different species are synchronized in
the sense that they are deterministic functions of the overlaps of the
whole system. This will describe the joint distribution of all overlaps
and allow us to obtain the lower bound in Section~\ref{Sec5label} in a
straightforward way using the Aizenman--Sims--Starr scheme. In the last
section, we will mention several interesting open questions.

\section{Guerra's replica symmetry breaking bound}\label{Sec2label}

Given $r\geq1$, let $(v_\alpha)_{\alpha\in\Natural^r}$ be the weights
of the Ruelle probability cascades \cite{Ruelle} corresponding to the
parameters (\ref{zetas}); see, for example, Section~2.3 in \cite{SKmodel} for the definition. For $\alpha, \beta\in\Natural^r$, we denote
%
\begin{equation}
\alpha\wedge\beta= \min \{0\leq\ell\leq r | \alpha_1= \beta
_1, \ldots, \alpha_{\ell} = \beta_{\ell},
\alpha_{\ell+1} \neq \beta _{\ell+1} \},
\end{equation}
where $\alpha\wedge\beta=r$ if $\alpha=\beta$. Since the sequences
defined in (\ref{Qs}) are nondecreasing, we can consider Gaussian
processes $C^s(\alpha)$ for $s\in\spe$ and $D(\alpha)$ both indexed by
$\alpha\in\Natural^r$ with the covariances
%
\begin{equation}
\e C^s(\alpha) C^s(\beta) = Q^s_{\alpha\wedge\beta}
\quad\mbox{and}\quad \e D(\alpha) D(\beta) = Q_{\alpha\wedge\beta}. \label{CD}
\end{equation}
These are the usual Gaussian fields that accompany the construction of
the Ruelle probability cascades; see, for example, Section~2.3 in \cite{SKmodel}. For each $s\in\spe$ and each $i\in I_s$, let $C_i(\alpha)$
be a copy of the process $C^s(\alpha)$, and suppose that all these
processes are independent of each other and of $D(\alpha)$. For $0\leq
x\leq1$, consider an interpolating Hamiltonian defined on $\Sigma
_N\times\Natural^r$ by
%
\begin{equation}
H_{N,x}(\sigma,\alpha) = \sqrt{x} H_N(\sigma) + \sqrt{1-x}
\smsp\sum_{i=1}^N \sigma_i
C_{i}(\alpha) +\sqrt{x} \sqrt{N} D(\alpha) \label{HNx}
\end{equation}
and the corresponding interpolating free energy
%
\begin{equation}
\varphi(x)=\frac{1}{N}\smsp\e\log\sum_{\sigma,\alpha}
v_{\alpha
} \exp H_{N, x}(\sigma,\alpha). \label{phix}
\end{equation}
Then it is easy to check the following.

\begin{lemma}
Under assumption (\ref{assumption}), the derivative of $\varphi(x)$ in
(\ref{phix}) satisfies $\varphi'(x)\leq0$.
\end{lemma}

\begin{pf} Let us denote by $\la \cdot\ra_x$ the average with
respect to the Gibbs measure $\Gamma_x(\sigma,\alpha)$ on $\Sigma
_N\times\Natural^r$ defined by
\[
\Gamma_x(\sigma,\alpha) \sim v_\alpha\exp
H_{N, x}(\sigma,\alpha).
\]
Then, obviously, for $0<x<1$,
\[
\varphi'(x) = \frac{1}{N}\e \biggl\la\frac{\partial H_{N, x}(\sigma
,\alpha)}{\partial x}
\biggr\ra_x.
\]
It is easy to check from the above definitions that
\begin{eqnarray*}
&& \frac{1}{N} \e\frac{\partial H_{N, x}(\sigma^1,\alpha^1)}{\partial
x} H_{N, x}\bigl(
\sigma^2,\alpha^2\bigr)
\\
&&\qquad = \frac{1}{2}\sum_{s,t\in\spe}
\Delta_{st}^2 \lambda_s \lambda_t
\bigl( R_s\bigl(\sigma^1,\sigma^2\bigr)
R_t\bigl(\sigma^1,\sigma^2\bigr)\\
&&\hspace*{106pt}{} -2
R_s\bigl(\sigma ^1,\sigma^2\bigr)
q^t_{\alpha^1\wedge\alpha^2} + q^s_{\alpha^1\wedge\alpha
^2}q^t_{\alpha
^1\wedge\alpha^2}
\bigr).
\end{eqnarray*}
In particular, this is zero when $(\sigma^1,\alpha^1)=(\sigma
^2,\alpha
^2)$ and, in general, can be rewritten as a quadratic form $(\Delta
^2(R-q),(R-q))/2$, where
\[
R= \bigl( \lambda_s R_s\bigl(\sigma^1,
\sigma^2\bigr) \bigr)_{s\in\spe},\qquad q= \bigl( \lambda_s
q^s_{\alpha^1\wedge\alpha^2} \bigr)_{s\in
\spe}.
\]
Notice that here we used the symmetry of the matrix $\Delta^2$.
Finally, usual Gaussian integration by parts then gives (see, e.g.,
Lemma~1.1 in \cite{SKmodel})
\[
\varphi'(x) = - \frac{1}{2}\e \bigl\la \bigl(
\Delta^2(R-q),(R-q) \bigr) \bigr\ra_x \leq0,
\]
where the last inequality follows from the assumption $\Delta^2\geq0$
in (\ref{assumption}).
\end{pf}

The lemma implies that $\varphi(1)\leq\varphi(0)$. It is easy to see that
\[
\varphi(0) = \log2 + \frac{1}{N}\smsp\e\log\sum
_{\alpha\in\Natural^r} v_{\alpha} \prod_{i\leq N}
\ch C_{i}(\alpha)
\]
and
\[
\varphi(1) = F_N + \frac{1}{N}\smsp\e\log\sum
_{\alpha\in\Natural^r} v_{\alpha} \exp\sqrt{N} D(\alpha).
\]
Now, standard properties of the Ruelle probability cascades imply that
(see, e.g., the proof of Lemma~3.1 in \cite{SKmodel}),
%
\begin{eqnarray}\label{simp1}
\frac{1}{N}\smsp\e\log\sum_{\alpha\in\Natural^r}
v_{\alpha} \prod_{i\leq N} \ch C_{i}(
\alpha) &=& \frac{1}{N}\smsp\sum_{1\leq i\leq N} \e\log\sum
_{\alpha\in
\Natural^r} v_{\alpha} \ch C_{i}(\alpha)
\nonumber
\\[-8pt]
\\[-8pt]
\nonumber
& =&
\sum_{s\in\spe} \lambda_s X^s_0
\end{eqnarray}
and
%
\begin{equation}
\frac{1}{N}\smsp\e\log\sum_{\alpha\in\Natural^r}
v_{\alpha} \exp\sqrt{N} D(\alpha) = \frac{1}{2} \sum
_{0\leq\ell\leq r-1} \zeta_\ell (Q_{\ell+1} -
Q_{\ell} ). \label{simp2}
\end{equation}
Recalling (\ref{Pzeta}), the inequality $\varphi(1)\leq\varphi(0)$ can
be written as $F_N\leq\PP(\zeta,q)$, which yields the upper bound in
(\ref{Parisi}).

\section{Multi-species Ghirlanda--Guerra identities} \label{Sec3label}

In order to prepare for the proof of the lower bound, we need to obtain
some strong coupling properties for the overlaps in different species,
which will be achieved in the next section using a multi-species
version of the Ghirlanda--Guerra identities that we will now prove. Let
us consider a countable dense subset $\WW$ of $[0,1]^{|\spe|}$. For a vector
%
\begin{equation}
w=(w_s)_{s\in\spe} \in\WW,
\end{equation}
let $s_i(w)=\sqrt{w_s}$ for $i\in I_s$ and $s\in\spe$, and consider the
following $p$-spin Hamiltonian,
%
\begin{equation}
h_{N,w,p}(\sigma) = \frac{1}{N^{p/2}} \sum_{1\leq i_1,\ldots,
i_p\leq
N}
g_{i_1,\ldots, i_p}^{w,p} \sigma_{i_1}s_{i_1}(w) \cdots
\sigma _{i_p}s_{i_p}(w), \label{hNwp}
\end{equation}
where $g_{i_1,\ldots, i_p}^{w,p}$ are i.i.d. standard Gaussian random
variables independent for all combinations of indices $p\geq1,w\in\WW
$ and $i_1,\ldots, i_p\in\{1,\ldots, N\}$. If we define
%
\begin{equation}
R_w \bigl(\sigma^1,\sigma^2 \bigr) = \sum
_{s\in\spe} \lambda_s w_s
R_s \bigl(\sigma^1,\sigma^2 \bigr),
\label{Rw}
\end{equation}
where $R_s(\sigma^1,\sigma^2)$ was defined in (\ref{overlap}), then it
is easy to check that the covariance of~(\ref{hNwp}) is
%
\begin{equation}
\e h_{N,w,p}\bigl(\sigma^1\bigr)h_{N,w,p}\bigl(
\sigma^2\bigr) = R_w \bigl(\sigma^1,
\sigma^2 \bigr)^p. \label{CovhNpw}
\end{equation}
Since the set $\WW$ is countable, we can consider some one-to-one
function \mbox{$j\dvtx \WW\to\Natural$}. Then we let $x_{w,p}$ for $p\geq1,w\in
\WW
$ be i.i.d. random variables uniform on the interval $[1,2]$ and define
a Hamiltonian
%
\begin{equation}
h_{N}(\sigma) = \sum_{w\in\WW} \sum
_{p\geq1} 2^{-j(w)-p} x_{w,p} h_{N,w,p}(
\sigma). \label{hNw}
\end{equation}
Note that, conditionally on $x=(x_{w,p})_{p\geq1,w\in\WW}$, this is a
Gaussian process and its variance is bounded by $4$. The Hamiltonian
$h_N(\sigma)$ will play a role of a perturbation Hamiltonian, which
means that, instead of $H_N(\sigma)$ in (\ref{SKH}), from now on we
will consider the perturbed Hamiltonian
%
\begin{equation}
H_N^{\mathrm{pert}}(\sigma) = H_N(\sigma) +
s_N h_N(\sigma), \label{Hpert}
\end{equation}
where $s_N=N^{\gamma}$ for any $1/4<\gamma<1/2$. First of all, it is
easy to see, using Jensen's inequality on each side, that
%
\begin{eqnarray}\label{compareF}
\frac{1}{N}\smsp\e\log\sum_{\sigma\in\Sigma_N} \exp
H_N(\sigma) & \leq& \frac{1}{N}\smsp\e\log\sum
_{\sigma\in\Sigma_N} \exp H_N^{\mathrm
{pert}}(\sigma)
\nonumber
\\[-8pt]
\\[-8pt]
\nonumber
& \leq& \frac{1}{N}\smsp\e\log\sum_{\sigma\in\Sigma_N} \exp
H_N(\sigma) + \frac{2s_N^2}{N},
\nonumber
\end{eqnarray}
and, since $\lim_{N\to\infty} N^{-1} s_N^2 = 0$, the perturbation term
does not affect the limit of the free energy. As in the
Sherrington--Kirkpatrick and mixed $p$-spin models, the purpose of
adding the perturbation term is to obtain the Ghirlanda--Guerra
identities for the Gibbs measure
%
\begin{equation}
G_N(\sigma) = \frac{\exp H_N^{\mathrm{pert}} (\sigma)}{Z_N} \qquad\mbox{where } Z_N =
\sum_{\sigma\in\Sigma_{N}} \exp H_N^{\mathrm
{pert}}(
\sigma), \label{GNpert}
\end{equation}
corresponding to the perturbed Hamiltonian (\ref{Hpert}). We will
denote the average with respect to $G_N^{\otimes\infty}$ by $\la
\cdot \ra$. Now, given $n\geq2$, let
\[
R^n = \bigl(R_s\bigl(\sigma^\ell,
\sigma^{\ell'}\bigr) \bigr)_{s\in\spe,
\ell,\ell
'\leq n}
\]
and consider an arbitrary bounded measurable function $f=f(R^n)$. For
$p\geq1$ and $w\in\WW$, let
%
\begin{eqnarray}\label{GGfinite}
\varDelta(f,n,w, p) &=&\biggl| {\e} \bigl\la f R_w \bigl(\sigma^1,
\sigma^{n+1} \bigr)^p\bigr \ra  - \frac{1}{n} {\e} \la f
\ra\smsp{\e} \bigl\la R_w \bigl(\sigma^1,\sigma^2
\bigr)^p \bigr\ra
\nonumber
\\[-8pt]
\\[-8pt]
\nonumber
&&\hspace*{89pt}{} - \frac{1}{n}\sum_{\ell=2}^{n}{\e}
\bigl\la f R_w \bigl(\sigma ^1,\sigma ^\ell
\bigr)^p \bigr\ra \biggr|,
\end{eqnarray}
where ${\e}$ denotes the expectation conditionally on the i.i.d.
uniform sequence $x=(x_{w,p})_{p\geq1,w\in\WW}$. If we denote by $\e
_x$ the expectation with respect to $x$ then the following holds.

\begin{theorem}\label{ThGG}
For any $n\geq2$ and any bounded measurable function $f=f(R^n)$,
%
\begin{equation}
\lim_{N\to\infty} \e_x \smsp\varDelta(f,n,w,p) = 0
\label{GGxlim}
\end{equation}
for all $p\geq1$ and $w\in\WW$.
\end{theorem}

\begin{pf}The proof is identical to the one of Theorem~3.2 in
\cite{SKmodel}. For a given $p\geq1$ and $w\in\WW$, equation (\ref
{GGxlim}) is obtained by utilizing the term $h_{N,w,p}(\sigma)$ in the
perturbation (\ref{hNw}).
\end{pf}

Theorem~\ref{ThGG} implies that we can choose a nonrandom sequence
$x^N= (x^N_{w,p})_{p\geq1,w\in\WW}$ changing with $N$ such that
%
\begin{equation}
\lim_{N\to\infty} \smsp\varDelta(f,n,w,p) = 0 \label{GGxlim2}
\end{equation}
for the Gibbs measure $G_N$ with the parameters $x$ in the perturbation
Hamiltonian~(\ref{hNw}) equal to $x^N$ rather than random. In fact, the
choice of $x^N$ will be made below in a special way to coordinate with
the Aizenman--Sim--Starr scheme. In this section, we will simply assume
that we have any such sequence $x^N$. Moreover, let us now consider any
subsequence $(N_k)_{k\geq1}$ along which the array
\[
\bigl(R_s\bigl(\sigma^\ell,\sigma^{\ell'}\bigr)
\bigr)_{s\in\spe,\ell
,\ell'\geq1}
\]
of the overlaps within species for infinitely many replicas $(\sigma
^\ell)_{\ell\geq1}$ converges in distribution under the measure $\e
G_N^{\otimes\infty}$. Again, later we will be interested in a special
choice of such subsequence. Let
%
\begin{equation}
\bigl(R^s_{\ell,\ell'}\bigr)_{s\in\spe,\ell,\ell'\geq1} \label{Rsfirst}
\end{equation}
be the array with the limiting distribution, and similarly to (\ref
{Rw}), define
%
\begin{equation}
R^w_{\ell,\ell'} = \sum_{s\in\spe}
\lambda_s w_s R^s_{\ell,\ell'}.
\label{Rwlim}
\end{equation}
Then equations (\ref{GGfinite}) and (\ref{GGxlim2}) imply that the
limiting array satisfies
%
\begin{equation}
\e f\bigl(R^n\bigr) \bigl(R^w_{1,n+1}
\bigr)^p = \frac{1}{n}\e f\bigl(R^n\bigr) \smsp\e
\bigl(R^w_{1,2} \bigr)^p + \frac{1}{n}\sum
_{\ell=2}^{n}\e f\bigl(R^n\bigr)
\bigl(R^w_{1,\ell} \bigr)^p, \label{GGwp}
\end{equation}
where, of course, now $R^n=(R^s_{\ell,\ell'})_{s\in\spe, \ell,\ell
'\leq
n}$. From this we will deduce the following multi-species form of the
Ghirlanda--Guerra identities for such limiting arrays. Let us consider
an array
%
\begin{equation}
Q_{\ell,\ell'} = \varphi \bigl( \bigl(R_{\ell,\ell'}^s
\bigr)_{s\in\spe
} \bigr) \label{Qphi}
\end{equation}
for any bounded measurable function $\varphi$ of the overlaps in
different species.

\begin{theorem}\label{ThGGms}
For any $n\geq2$ and any bounded measurable function $f=f(R^n)$,
%
\begin{equation}
\e f\bigl(R^n\bigr) Q_{1,n+1} = \frac{1}{n}\e f
\bigl(R^n\bigr) \smsp\e Q_{1,2} + \frac{1}{n}\sum
_{\ell=2}^{n}\e f\bigl(R^n\bigr)
Q_{1,\ell}. \label{GGms}
\end{equation}
\end{theorem}

\begin{pf} Since equation (\ref{GGwp}) holds for all $w\in\WW$,
both sides are continuous in $w$, and $\WW$ is dense in $[0,1]^{|\spe
|}$, equation (\ref{GGwp}) holds for all $w\in[0,1]^{|\spe|}$. Take
any integers $p_s\geq0$ for $s\in\spe$, and let $p=\sum_{s\in\spe}
p_s$. If we recall the definition of $R^w_{\ell,\ell'}$ in (\ref{Rwlim}),
\[
\frac{\partial^p}{\prod_{s\in\spe} \partial w_s^{p_s}} \bigl(R^w_{\ell
,\ell'} \bigr)^p =
p! \prod_{s\in\spe} \bigl(\lambda_s
R^s_{\ell,\ell'}\bigr)^{p_s}.
\]
Computing this partial derivative on both sides of (\ref{GGwp}) implies
%
\begin{eqnarray}
\e f\bigl(R^n\bigr) \prod_{s\in\spe}
\bigl(R^s_{1,n+1}\bigr)^{p_s}& = & \frac{1}{n}\e
f\bigl(R^n\bigr) \smsp\e\prod_{s\in\spe}
\bigl(R^s_{1,2}\bigr)^{p_s}
\nonumber
\\[-8pt]
\\[-8pt]
\nonumber
&&{} + \frac{1}{n}\sum_{\ell=2}^{n}\e f
\bigl(R^n\bigr) \prod_{s\in\spe}
\bigl(R^s_{1,\ell}\bigr)^{p_s}.
\nonumber
\end{eqnarray}
Approximating continuous functions by polynomials, this implies (\ref
{GGms}) for continuous functions $\varphi$ in (\ref{Qphi}), and the
general case follows.
\end{pf}

\begin{remark*} In particular, Theorem~\ref{ThGGms} implies that the
array $(Q_{\ell,\ell'})_{\ell,\ell'\geq1}$ itself satisfies the usual
Ghirlanda--Guerra identities,
%
\begin{equation}\quad
\e f\bigl(Q^n\bigr) \psi(Q_{1,n+1}) = \frac{1}{n}\e f
\bigl(Q^n\bigr) \smsp\e\psi(Q_{1,2}) + \frac{1}{n}\sum
_{\ell=2}^{n}\e f\bigl(Q^n\bigr)
\psi(Q_{1,\ell}), \label{GGreg}
\end{equation}
for any bounded measurable function $\psi$ and $f=f(Q^n)$, where $Q^n
=\break  (Q_{\ell,\ell'})_{\ell,\ell'\leq n}$. In the case when the array $Q$
is also nonnegative definite, the main result in \cite{PUltra} will
allow us to use the full force of the Ghirlanda--Guerra identities and,
in particular, will imply that such arrays are ultrametric and can be
generated by the Ruelle probability cascades;
see Section~2.4 in \cite{SKmodel}.
\end{remark*}

\section{Synchronizing the species}\label{Sec4label}

Now, let us consider any limiting distribution as in (\ref{Rsfirst}),
and let us notice that the overlap
\[
R\bigl(\sigma^\ell,\sigma^{\ell'}\bigr) = \frac{1}{N}
\sum_{i=1}^N \sigma ^\ell_i
\sigma^{\ell'}_i = \sum_{s\in\spe}
\lambda_s R_s\bigl( \sigma^\ell,
\sigma^{\ell'}\bigr)
\]
of two configurations over the whole system in the limit will become
%
\begin{equation}
R_{\ell,\ell'} = \sum_{s\in\spe} \lambda_s
R^s_{\ell,\ell'}. \label{Aoverlap}
\end{equation}
In this section, we will prove the main result that will allow us to
characterize the limits that will arise in the Aizenman--Sims--Starr scheme.

\begin{theorem}\label{ThSynch}
For any array (\ref{Rsfirst}) that satisfies (\ref{GGms}), there exist
nondecreasing $(1/\lambda_s)$-Lipschitz functions $L_s\dvtx [0,1]\to
[0,1]$ such that $R^s_{\ell,\ell'} = L_s(R_{\ell,\ell'})$ almost surely
for all $s\in\spe$ and all $\ell,\ell'\geq1$.
\end{theorem}

The reason we can consider the domain and range of $L_s$ to be $[0,1]$
is because each array $R^s$ is nonnegative definite and satisfies the
Ghirlanda--Guerra identities~(\ref{GGreg}), and therefore, its entries
are nonnegative by Talagrand's positivity principle (Theorem~2.16 in
\cite{SKmodel}). Theorem~\ref{ThSynch} implies that the joint
distribution of the overlap arrays for all species will be determined
trivially by the overlap array $(R_{\ell,\ell'})_{\ell,\ell'\geq1}$.
On the other hand, the Ghirlanda--Guerra identities imply that this
array can be generated using the Ruelle probability cascades, which
will be used in Section~\ref{Sec5label}. We begin with the following
observation.

\begin{lemma}\label{LemSynch}
If $R^s_{\ell,\ell'}> R^s_{\ell,\ell''}$ for some $s\in\spe$, then
$R^t_{\ell,\ell'}\geq R^t_{\ell,\ell''}$ for all $t\in\spe$.
\end{lemma}

\begin{pf}
By Theorem~\ref{ThGGms}, for any $s,t\in\spe$, the arrays
\[
\bigl(R^s_{\ell,\ell'} \bigr)_{\ell,\ell'\geq1},\qquad
\bigl(R^t_{\ell
,\ell
'} \bigr)_{\ell,\ell'\geq1}\quad \mbox{and}\quad
\bigl(R^s_{\ell,\ell'}+R^t_{\ell,\ell'}
\bigr)_{\ell,\ell'\geq1}
\]
satisfy the Ghirlanda--Guerra identities. Since all these arrays are
nonnegative definite, the main result in \cite{PUltra} (or Theorem~2.14
in \cite{SKmodel}) implies that these arrays are ultrametric, that is,
%
\begin{equation}
R^s_{\ell',\ell''} \geq\min \bigl(R^s_{\ell,\ell'},
R^s_{\ell
,\ell''} \bigr) \label{ultram}
\end{equation}
for any different $\ell,\ell',\ell'' \geq1$ and, similarly, for the
other two arrays. In other words, given three replica indices, the
smallest two overlaps are equal. Suppose now that $R^s_{\ell,\ell'}>
R^s_{\ell,\ell''}$ but $R^t_{\ell,\ell'}< R^t_{\ell,\ell''}$. By
ultrameticity of the first two arrays,
\[
R^s_{\ell,\ell'}> R^s_{\ell,\ell''} =
R^s_{\ell',\ell''}\quad \mbox{and}\quad R^t_{\ell',\ell''} =
R^t_{\ell,\ell'}< R^t_{\ell,\ell''}.
\]
However, this implies that
\[
R^s_{\ell',\ell''} + R^t_{\ell',\ell''}< \min \bigl(
R^s_{\ell
,\ell
'}+R^t_{\ell,\ell'},
R^s_{\ell,\ell''} + R^t_{\ell,\ell''} \bigr),
\]
violating ultrametricity of the third array.
\end{pf}

Let us state one obvious corollary of the above lemma.

\begin{corollary}\label{CorSynch}
The following statements hold:
\begin{longlist}[(a)]
\item[(a)] If $R_{\ell,\ell'}> R_{\ell,\ell''}$, then $R^s_{\ell
,\ell
'}\geq R^s_{\ell,\ell''}$ for all $s\in\spe$.

\item[(b)] If $R^s_{\ell,\ell'}> R^s_{\ell,\ell''}$ for some $s\in
\spe
$, then $R_{\ell,\ell'}> R_{\ell,\ell''}$.
\end{longlist}
\end{corollary}

This already gives some indication that the overlaps in different
species will be synchronized. However, keeping in mind the ultrametric
tree structure of the Ruelle probability cascades that generate them,
we need to show that the entire clusters are synchronized and the
corresponding cascades are completely coupled. To prove this, for $q\in
[0,1]$ and $s\in\spe$, we will consider the array
%
\begin{equation}
R^{s,q}_{\ell,\ell'} = \Ind(R_{\ell,\ell'}\geq q)
\bigl(R^s_{\ell,\ell
'} +1\bigr).
\end{equation}
First of all, we add $+1$ to the overlap $R^s_{\ell,\ell'}$ to ensure
that the only way the right-hand side can be equal to zero is when
$R_{\ell,\ell'}< q$ and not, for example, when $R^s_{\ell,\ell
'}=0$. As
in (\ref{ultram}), by Theorem~\ref{ThGGms}, the array $(R_{\ell,\ell
'})_{\ell,\ell'\geq1}$ is ultrametric, which implies that the array
$(\Ind(R_{\ell,\ell'}\geq q))_{\ell,\ell'\geq1}$ is nonnegative
definite, as it consists of blocks on the diagonal with all entries
equal to one. Therefore, the array
\[
R^{s,q} = \bigl(R^{s,q}_{\ell,\ell'} \bigr)_{\ell,\ell'\geq1}
\]
is nonnegative definite as the Hadamard product of two such arrays. By
Theorem~\ref{ThGGms}, the array $R^{s,q}$ also satisfies the
Ghirlanda--Guerra identities, so all the consequences of the
Ghirlanda--Guerra identities for nonnegative definite arrays described,
for example, in Section~2.4 in \cite{SKmodel}, hold in this case. One
such consequence is the following. Let
%
\begin{equation}
\mu= \LL(R_{1,2}) \quad\mbox{and}\quad \mu^{s,q} = \LL
\bigl(R^{s,q}_{1,2}\bigr)
\end{equation}
be the distributions of one entry of the arrays $R$ and $R^{s,q}$
correspondingly. Lemma~2.7 in \cite{SKmodel} implies the following
consequence of the Ghirlanda--Guerra identities, which was first
observed in \cite{PT}.

\begin{lemma}\label{LemDense}
For any $s\in\spe, \ell\geq1$ and $q\in[0,1]$, with probability one,
the set
%
\begin{equation}
A^{s}_\ell(q) = \bigl\{R^{s,q}_{\ell,\ell'} |
\ell'\neq\ell \bigr\} = \bigl\{ \Ind(R_{\ell,\ell'}\geq q)
\bigl(R^s_{\ell,\ell'} +1\bigr) | \ell '\neq\ell \bigr
\} \label{Asq}
\end{equation}
is a dense subset of the support of $\mu^{s,q}$.
\end{lemma}

This will be the key to the proof of Theorem~\ref{ThSynch}. Now, for
any $q\in[0,1]$, let us define
%
\begin{equation}
\ell_s(q) = \inf \bigl\{x\geq1 | x\in\operatorname{supp} \mu^{s,q}
\bigr\}-1. \label{ellsq}
\end{equation}
Equivalently, one could take the infimum over $x>0$, because
$R^{s,q}_{\ell,\ell'}>0$ if and only if $R^{s,q}_{\ell,\ell'}\geq1$.
To understand the meaning of this definition, let us notice that,
whenever the set $A^{s}_\ell(q)$ in (\ref{Asq}) is dense in the support
of $\mu^{s,q}$ (which happens with probability one for a given $q$),
%
\begin{equation}
\ell_s(q) = \inf \bigl\{R^s_{\ell,\ell'} |
\ell'\neq\ell, R_{\ell
,\ell'} \geq q \bigr\}, \label{ellsq2}
\end{equation}
so $\ell_s(q)$ is just the smallest value that $R^s_{\ell,\ell'}$ can
take whenever $R_{\ell, \ell'}\geq q$. This alternative definition,
obviously, implies the following.

\begin{lemma}
For any $s\in\spe$, the function $\ell_s(q)$ in (\ref{ellsq}) is
nondecreasing in $q$.
\end{lemma}

To obtain the functions $L_s$ in Theorem~\ref{ThSynch}, we will first
need to regularize $\ell_s(q)$ as follows:
%
\begin{equation}
L_s(q) = \lim_{x\uparrow q} \ell_s(x)
\end{equation}
for $q>0$ and $L_s(0) = \ell_s(0)$. Theorem~\ref{ThSynch} will be now
proved in two steps. First, we will show that $R^s_{\ell,\ell'} =
L_s(R_{\ell,\ell'})$ almost surely. Second, we will show that $L_s$ is
$(1/\lambda_s)$-Lipschitz on the support of the distribution $\mu$ of
$R_{1,2}$. Then we can redefine $L_s$ outside of the support to be
$(1/\lambda_s)$-Lipschitz extension which, obviously, does not change
the first claim, $R^s_{\ell,\ell'} = L_s(R_{\ell,\ell'})$, since
$R_{\ell,\ell'}$ belongs to the support of $\mu$ almost surely.

\begin{pf*}{Proof of Theorem~\ref{ThSynch}}
\emph{Step} 1. We will use
that the claim in Lemma~\ref{LemDense} holds with probability one
simultaneously for all $q\in\q\cap[0,1]$. Let us fix some indices
$\ell\neq \ell'$. If $\mu(\{0\})=0$, then all $R_{\ell,\ell'} > 0$
almost surely. If $\mu(\{0\})>0$ and $R_{\ell,\ell'} = 0$, then we must
have $R^s_{\ell,\ell'}=0$ for all $s\in\spe$, and definition (\ref
{ellsq}) implies that $\ell_s(0)=0$. In this case,
\[
R^s_{\ell,\ell'} = \ell_s(R_{\ell,\ell'})
=L_s(R_{\ell,\ell'}).
\]
Let us now consider the case when $R_{\ell,\ell'}>0$. First of all, for
any $x< R_{\ell,\ell'}$ we must have that $\ell_s(x)\leq R^{s}_{\ell
,\ell'}$, because the function $\ell_s(x)$ is nondecreasing and, for
any rational $q\leq R_{\ell,\ell'}$, (\ref{ellsq2}) implies that
$\ell
_s(q)\leq R^{s}_{\ell,\ell'}$. Next, consider arbitrary $\eps>0$, and
consider any rational $q$ such that
%
\begin{equation}
q<R_{\ell,\ell'}\leq q+\eps. \label{qeps}
\end{equation}
Consider two possibilities. First, suppose that $R^s_{\ell,\ell'} =
\ell
_s(q)$. Since for $q\leq x < R_{\ell,\ell'}$ we showed that
\[
\ell_s(x) \leq R^s_{\ell,\ell'} =
\ell_s(q) \leq\ell_s(x)
\]
[so $\ell_s(x) = \ell_s(q)$ for such $x$], we get the desired claim,
\[
R_{\ell,\ell'}^s = \ell_s(q) = \lim
_{x\uparrow R_{\ell,\ell'}} \ell _s(x) = L_s(R_{\ell,\ell'}).
\]
Second, suppose that $\ell_s(q)< R^s_{\ell,\ell'}$. By (\ref{ellsq2}),
we can find a sequence $(\ell_n)$ such that $R_{\ell,\ell_n}\geq q$ and
$R^{s}_{\ell,\ell_n} \downarrow\ell_s(q)$. Since we assumed that
$\ell
_s(q)< R^s_{\ell,\ell'}$, for large enough $n$ we must have
$R^{s}_{\ell
,\ell_n}< R^{s}_{\ell,\ell'}$ and, by Corollary~\ref{CorSynch}, we get
$R_{\ell,\ell_n}< R_{\ell,\ell'}$ and $R^{t}_{\ell,\ell_n}\leq
R^{t}_{\ell,\ell'}$ for all $t\in\spe$. Therefore,
%
\begin{eqnarray} \label{Lip}
0&\leq&\lambda_s \bigl(R^{s}_{\ell,\ell'} -
R^{s}_{\ell,\ell_n} \bigr) \leq\sum_{t\in\spe}
\lambda_t \bigl(R^{t}_{\ell,\ell'} -
R^{t}_{\ell
,\ell_n} \bigr)
\nonumber
\\[-8pt]
\\[-8pt]
\nonumber
&=& R_{\ell,\ell'} - R_{\ell,\ell_n}
\leq q+\eps-q = \eps.
\end{eqnarray}
Using that $R^{s}_{\ell,\ell_n} \downarrow\ell_s(q)$ implies that
$\ell
_s(q) \leq R^{s}_{\ell,\ell'} \leq\ell_s(q) + \eps\lambda_s^{-1}$.
Finally, letting $q\uparrow R_{\ell,\ell'}$ and $\eps\downarrow0$ in
such a way that (\ref{qeps}) holds, again, implies the desired claim
\[
R_{\ell,\ell'}^s = \lim_{q\uparrow R_{\ell,\ell'}}
\ell_s(q) = L_s(R_{\ell,\ell'}).
\]
\emph{Step} 2. Let us now show that $L_s$ is $(1/\lambda_s)$-Lipschitz
on the support of the distribution $\mu$ of $R_{1,2}$. Take $q_1< q_2$
in the support of $\mu$. Let $q_2' = q_2-\eps_2$ for some small $\eps
_2>0$ such that $q_2'>0$, and let $q_1'=\max(q_1-\eps_1, 0)$ for some
small $\eps_1>0$. Let us also make sure that $q_1'$ and $q_2'$ are
rational. By (\ref{ellsq2}), given $\eps>0$, we can find indices
$\ell
_j$ for $j=1,2$ such that
%
\begin{equation}
R_{\ell,\ell_j}\geq q_j'\quad \mbox{and}\quad
\ell_s\bigl(q_j'\bigr) \leq
R^s_{\ell,\ell_j}\leq\ell_s\bigl(q_j'
\bigr)+\eps. \label{temp}
\end{equation}
Similarly to Lemma~\ref{LemDense}, Lemma~2.7 in \cite{SKmodel} implies
that the set $\{R_{\ell,\ell'} | \ell'\neq\ell\}$ is a dense subset
of the support of $\mu= \LL(R_{1,2})$ with probability one and, since
we chose $q_1$ and $q_2$ in the support of $\mu$, we can find other
indices $\ell_j'$ for $j=1,2$ such that
\[
q_j' \leq R_{\ell,\ell_j'} \leq q_j+
\eps.
\]
If the index $\ell_j$ already satisfies this condition, we simply take
$\ell_j' = \ell_j$. Otherwise, because of the first inequality in
(\ref
{temp}), we must have $R_{\ell,\ell_j'} < R_{\ell,\ell_j}$ and, by
(\ref
{ellsq2}), Corollary~\ref{CorSynch} and the second inequality in (\ref{temp}),
\[
\ell_s\bigl(q_j'\bigr) \leq
R^s_{\ell,\ell_j'} \leq R^s_{\ell,\ell_j} \leq \ell
_s\bigl(q_j'\bigr)+\eps.
\]
In both cases, we have
\[
q_j' \leq R_{\ell,\ell_j'} \leq q_j+
\eps \quad\mbox{and}\quad \ell_s\bigl(q_j'\bigr)
\leq R^s_{\ell,\ell_j'} \leq\ell_s\bigl(q_j'
\bigr)+\eps.
\]
Since $q_1<q_2$, by taking $\eps>0$ small enough, we can assume that
$R_{\ell,\ell_1'}<R_{\ell,\ell_2'}$. Then, as in (\ref{Lip}),
\[
\lambda_s \bigl( R^{s}_{\ell,\ell_2'} -
R^{s}_{\ell,\ell_1'} \bigr) \leq R_{\ell,\ell_2'} -
R_{\ell,\ell_1'}.
\]
Combining all the inequalities, we showed that
\[
\lambda_s \bigl( \ell_s\bigl(q_2'
\bigr) - \ell_s\bigl(q_1'\bigr)-\eps \bigr)
\leq q_2 +\eps-q_1'.
\]
Letting $\eps, \eps_1, \eps_2\downarrow0$ implies $\lambda_s (L_s(q_2)
- L_s(q_1)) \leq q_2 - q_1$, which proves that $L_s$ is $(1/\lambda
_s)$-Lipschitz on the support of $\mu$. As we mentioned above,
$(1/\lambda_s)$-Lipschitz extension of $L_s$ outside of the support
does not affect the fact that $R^s_{\ell,\ell'} = L_s(R_{\ell,\ell'})$
almost surely.
\end{pf*}

\section{Lower bound via the Aizenman--Sims--Starr scheme}\label{Sec5label}

Given the main result in the previous section, the arguments of this
section will be a standard exercise. To a reader familiar with the
corresponding arguments in the setting of the classical SK model (e.g.,
Sections~3.5 and 3.6 in \cite{SKmodel}) these arguments will be
completely obvious. Otherwise, we recommend to study them first in the
easier case of the SK model.

It is clear that small modifications of the vector $(\lambda_s)_{s\in
\spe}$ result in small changes both of the free energy for large $N$
and the Parisi formula (\ref{Parisi}), so without loss of generality,
we can assume that all $\lambda_s$ are rational and can be written as
%
\begin{equation}
\lambda_s = \frac{k_s}{k}.
\end{equation}
In the proof of the lower bound, we will use an obvious fact that
%
\begin{equation}
\liminf_{N\to\infty} F_N \geq\frac{1}{k}\liminf
_{n\to\infty} (\e \log Z_{nk+k} - \e\log Z_{nk} ).
\label{FNAN}
\end{equation}
Let us consider the right-hand side for a fixed $N=nk$, and in addition
to partition~(\ref{species}), let us consider a partition of $k$ new
coordinates
%
\begin{equation}
I^+ = \{N+1,\ldots, N+k\} = \bigcup_{s\in\spe}
I_s^+ \label{speciescavity}
\end{equation}
into different species, so that $|I_s^+| = k_s$. Let us compare the
partition functions $Z_{N}$ and $Z_{N+k}$. If we denote $\rho= (\sigma
,\eps)\in\Sigma_{N+k}$ for $\sigma\in\Sigma_N$ and $\eps
\in\Sigma
_k$, then we can write
%
\begin{equation}
H_{N+k}(\rho) = H_N'(\sigma) + \sum
_{i\in I^+}\eps_i z_{N,i}(\sigma ) + r(
\eps), \label{decomp1}
\end{equation}
where
%
\begin{eqnarray}
H_N'(\sigma)& = &\frac{1}{\sqrt{N+k}} \sum
_{i,j =1}^N g_{ij}\sigma_i
\sigma_j, \label{commonH}
\\
z_{N,i}(\sigma) &=& \frac{1}{\sqrt{N+k}} \sum
_{j=1}^N (g_{ij} + g_{ji} )
\sigma_j
\end{eqnarray}
and
%
\begin{equation}
r(\eps)= \frac{1}{\sqrt{N+k}} \sum_{i,j\in I^+}
g_{ij}\eps_i \eps_j.
\end{equation}
On the other hand, the Gaussian process $H_N(\sigma)$ on $\Sigma_N$ can
be decomposed into a sum of two independent Gaussian processes
%
\begin{equation}
H_N(\sigma) \stackrel{d} {=} H_N'(\sigma)
+ y_N(\sigma), \label{commonH2}
\end{equation}
where
%
\begin{equation}
y_N(\sigma) = \frac{\sqrt{k}}{\sqrt{N(N+k)}} \sum_{i,j =1}^N
g_{ij}'\sigma_i \sigma_j
\end{equation}
and $(g_{ij}')$ are independent copies of the Gaussian random variables
$(g_{ij})$. Using that the term $r(\eps)$ is of a small order, we can write
%
\begin{equation}
\e\log Z_{N+k} = \e\log\sum_{\sigma\in\Sigma_N} \prod
_{i\in I^+} 2\smsp\ch \bigl(z_{N,i}(\sigma)
\bigr) \exp H_{N}'(\sigma) + o(1) \label{ZN1}
\end{equation}
and, using equation (\ref{commonH2}),
%
\begin{equation}
\e\log Z_{N} = \e\log\sum_{\sigma\in\Sigma_N} \exp
\bigl(y_N(\sigma) \bigr) \exp H_{N}'(\sigma).
\label{ZN}
\end{equation}
Finally, if we consider the Gibbs measure on $\Sigma_N$
corresponding to the Hamiltonian $H_N'(\sigma)$ in (\ref{commonH}),
%
\begin{equation}
G_N'(\sigma) = \frac{\exp H_N'(\sigma)}{Z_N'}\qquad \mbox{where }
Z_N' = \sum_{\sigma\in\Sigma_{N}} \exp
H_N'(\sigma), \label{MeasureGNprime}
\end{equation}
then combining (\ref{ZN1}), (\ref{ZN}) we can replace the right-hand
side of (\ref{FNAN}) by
%
\begin{eqnarray}\label{AS2repr}
&&\frac{1}{k}\liminf_{n\to\infty} \biggl(  \e\log\sum
_{\sigma\in\Sigma_{N}} \prod_{i\in I^+} 2\smsp\ch
\bigl(z_{N,i}(\sigma) \bigr) G_N'(\sigma)
\nonumber
\\[-8pt]
\\[-8pt]
\nonumber
&&\hspace*{58pt}{}- \e\log\sum_{\sigma\in\Sigma_{N}} \exp \bigl(y_N(
\sigma) \bigr) G_N'(\sigma) \biggr).
\end{eqnarray}
This is the analogue of the Aizenman--Sims--Starr representation in
\cite{AS2}; see Section~3.5 in \cite{SKmodel}. From the construction it
is clear that the Gaussian processes $z_{N,i}(\sigma)$ for $i\in I^+$
and $y_N(\sigma)$ are independent of each other and the randomness of
the measure $G_N'$. For $s\in\spe$ and $i\in I_s^+$,
%
\begin{eqnarray}\label{Covz}
\e z_{N,i}\bigl(\sigma^1\bigr) z_{N,i}\bigl(
\sigma^2\bigr) &= &\frac{1}{N+k}\sum_{t\in\spe}
\sum_{j\in I_t} 2 \Delta_{st}^2
\sigma^1_j \sigma^2_j\nonumber
\\
&=& \frac{N}{N+k}\sum_{t\in\spe} 2
\Delta_{st}^2 \lambda_t R_t\bigl(
\sigma^1, \sigma^2\bigr)
\\
&=& 2\sum_{t\in\spe} \Delta_{st}^2
\lambda_t R_t\bigl(\sigma^1,
\sigma^2\bigr) + O\bigl(N^{-1}\bigr)
\nonumber
\end{eqnarray}
and, similarly to the computation of the covariance in (\ref{Cov}),
%
\begin{eqnarray}\label{Covy}
\e y_{N}\bigl(\sigma^1\bigr) y_{N}\bigl(
\sigma^2\bigr) &=& \frac{kN^2}{N(N+k)}\sum_{s,t\in\spe}
\Delta_{st}^2 \lambda_s \lambda_t
R_s\bigl(\sigma^1, \sigma^2\bigr)
R_t\bigl(\sigma^1, \sigma^2\bigr)
\nonumber
\\[-8pt]
\\[-8pt]
\nonumber
&=& k \sum_{s,t\in\spe} \Delta_{st}^2
\lambda_s \lambda_t R_s\bigl(
\sigma^1, \sigma^2\bigr) R_t\bigl(
\sigma^1, \sigma^2\bigr) + O\bigl(N^{-1}\bigr).
\nonumber
\end{eqnarray}
Notice how these expressions resemble the definition in (\ref{Qs}). Of
course, one can ignore the lower error terms $O(N^{-1})$ from now on.

The same computation can be carried out just as easily in the case when
the free energy $F_N$ in (\ref{FNAN}) corresponds to the perturbed
Hamiltonian $H_N^{\mathrm{pert}}(\sigma)$ in (\ref{Hpert}) instead of
the original Hamiltonian $H_N(\sigma)$. Moreover, since the
perturbation term $s_N h_N(\sigma)$ in (\ref{Hpert}) is of a smaller
order, one can show that the perturbation term $s_{N+k} h_{N+k}(\rho)$
in the partition function $Z_{N+k}$ can simply be replaced by the one
in $Z_{N}$, $s_N h_N(\sigma)$. This is standard and is explained, for
example, in Section~3.5 in~\cite{SKmodel}. In this case, we obtain the
representation (\ref{AS2repr}) with the Gibbs measure $G_N'$ in (\ref
{MeasureGNprime}) corresponding to the perturbed Hamiltonian
\[
H_N'(\sigma) + s_N h_N(
\sigma).
\]
Also, in this case the expectation $\e$ in (\ref{AS2repr}) includes the
average $\e_x$ in the uniform random variables $x=(x_{w,p})$ in the
definition of the perturbation Hamiltonian (\ref{hNw}).

The proof of Theorem~\ref{ThGG} applies verbatim to the measure $G_N'$,
and right below Theorem~\ref{ThGG} we mentioned that one can choose a
nonrandom sequence $x^N=(x^N_{w,p})_{p\geq1,w\in\WW}$ changing with
$N$ such that (\ref{GGxlim2}) holds for the Gibbs measure $G_N'$ with
the parameters $x$ in the perturbation Hamiltonian (\ref{hNw}) equal to
$x^N$ rather than random. By Lemma~3.3 in \cite{SKmodel}, one can
choose this sequence $x^N$ in such a way that the lower limit in (\ref
{AS2repr}) is not affected by fixing $x=x^N$ instead of averaging in~$x$.
To finish the proof, we will use Theorem~1.3 in \cite{SKmodel} (a
trivial modification of) which implies that
%
\begin{equation}
\e\log\sum_{\sigma\in\Sigma_{N}} \prod
_{i\in I^+} 2\smsp\ch \bigl(z_{N,i}(\sigma) \bigr)
G_N'(\sigma) - \e\log\sum_{\sigma\in\Sigma_{N}}
\exp \bigl(y_N(\sigma) \bigr) G_N'(\sigma)
\label{comp1}
\end{equation}
is a continuous functional of the distribution of the array
%
\begin{equation}
\bigl(R_s\bigl(\sigma^\ell,\sigma^{\ell'}\bigr)
\bigr)_{s\in\spe,\ell
,\ell'\geq1} \label{Rsend}
\end{equation}
under the measure $\e G_N^{\prime\otimes\infty}$. Passing to a
subsequence, if necessary, we can assume that this array converges in
distribution to some array $(R^s_{\ell,\ell'})_{s\in\spe,\ell,\ell
'\geq
1}$ that, by construction, satisfies Theorem~\ref{ThGGms}. In
particular, by Theorem~\ref{ThSynch},
%
\begin{equation}
R^s_{\ell,\ell'} = L_s(R_{\ell,\ell'})
\label{RLS}
\end{equation}
for some nondecreasing $(1/\lambda_s)$-Lipschitz functions $L_s$,
where $R_{\ell,\ell'}$ is the overlap of the whole system in (\ref
{Aoverlap}).

Let us consider sequence (\ref{zetas}) and a sequence
%
\begin{equation}
0=q_0< q_1 < \cdots< q_{r-1}< q_r
=1 \label{q}
\end{equation}
such that the distribution $\zeta$ on $[0,1]$ defined by
%
\begin{equation}
\zeta \bigl( \{q_\ell \} \bigr) = \zeta_{\ell} - \zeta
_{\ell-1} \qquad\mbox{for } \ell=0,\ldots, r \label{zetafop}
\end{equation}
is close to the distribution $\LL(R_{1,2})$ of one element of the array
$(R_{\ell,\ell'})_{\ell,\ell'\geq1}$ in some metric that metrizes weak
convergence of distributions on $[0,1]$. As in Section~\ref{Sec2label},
let $(v_\alpha)_{\alpha\in\Natural^r}$ be the weights of the Ruelle
probability cascades corresponding to the parameters (\ref{zetas}). Let
$(\alpha^\ell)_{\ell\geq1}$ be an i.i.d. sample from $\Natural^r$
according to these weights and, using sequence (\ref{q}), define
%
\begin{equation}
Q_{\ell,\ell'} = q_{\alpha^\ell\wedge\alpha^{\ell'}}.
\end{equation}
Since from Theorem~\ref{ThGGms} it is clear that the overlap array
$(R_{\ell,\ell'})_{\ell,\ell'\geq1}$ satisfies the Ghirlanda--Guerra
identities, Theorems 2.13 and 2.17 in \cite{SKmodel} imply that its
distribution will be close to the distribution of the array $(Q_{\ell
,\ell'})_{\ell,\ell'\geq1}$. If for each $s\in\spe$ we define the
sequence in (\ref{qs}) by
%
\begin{equation}
q^s_\ell= L_s(q_\ell)\qquad \mbox{for } 0
\leq l\leq r, \label{zetafops}
\end{equation}
and let
%
\begin{equation}
Q^s_{\ell,\ell'} = L_s(Q_{\ell,\ell'}) =
q^s_{\alpha^\ell\wedge
\alpha
^{\ell'}}, \label{QLS}
\end{equation}
equation (\ref{RLS}) implies that the entire array $(Q^s_{\ell,\ell
'})_{s,\in\spe, \ell,\ell'\geq1}$ will be close in distribution to the
array $(R^s_{\ell,\ell'})_{s\in\spe, \ell,\ell'\geq1}$.

Let us now consider Gaussian processes $C^s(\alpha)$ for $s\in\spe$ and
$D(\alpha)$ indexed by $\alpha\in\Natural^r$ as in Section~\ref{Sec2label}. For each $s\in\spe$ and each $i\in I^+_s$, let
$C_i(\alpha
)$ be a copy of the process $C^s(\alpha)$, and suppose that all these
processes are independent of each other and of $D(\alpha)$. Similarly
to (\ref{comp1}), consider
%
\begin{equation}
\e\log\sum_{\alpha\in\Natural^r} \prod_{i\in I^+}
2\smsp\ch \bigl(C_{i}(\alpha) \bigr) v_\alpha - \e\log\sum
_{\alpha\in\Natural^r} \exp \bigl(\sqrt{k}D(\alpha ) \bigr)
v_\alpha. \label{comp2}
\end{equation}
By (\ref{Qs}), (\ref{CD}) and (\ref{QLS}), the covariances of these
Gaussian processes can be written as
%
\begin{equation}
\e C_i\bigl(\alpha^1\bigr) C_i\bigl(
\alpha^2\bigr) = 2\sum_{t\in\spe}
\Delta_{st}^2 \lambda_t q^s_{\alpha^1\wedge
\alpha^2}
= 2\sum_{t\in\spe} \Delta_{st}^2
\lambda_t Q^s_{1,2} \label{CD1}
\end{equation}
for $s\in\spe$ and $i\in I_s^+$, and
%
\begin{eqnarray}\label{CD2}
\e\sqrt{k}D\bigl(\alpha^1\bigr) \sqrt{k}D\bigl(\alpha^2
\bigr) &=& k \sum_{s,t\in\spe} \Delta_{st}^2
\lambda_s \lambda_t q^s_{\alpha
^1\wedge\alpha^2}q^t_{\alpha^1\wedge\alpha^2}
\nonumber
\\[-8pt]
\\[-8pt]
\nonumber
&=& k \sum_{s,t\in\spe} \Delta_{st}^2
\lambda_s \lambda_t Q^s_{1,2}
Q^t_{1,2}.
\nonumber
\end{eqnarray}
If we compare the covariances in (\ref{Covz}) and (\ref{Covy}) with
(\ref{CD1}) and (\ref{CD2}), Theorem~1.3 in \cite{SKmodel} implies that
(\ref{comp2}) is the same continuous functional of the distribution of
the array
%
\begin{equation}
\bigl(Q^s_{\ell,\ell'} \bigr)_{s\in\spe,\ell,\ell'\geq1} \label{Qsend},
\end{equation}
as (\ref{comp1}) is of the array (\ref{Rsend}). Since both arrays, by
construction, approximate in distribution the array $(R^s_{\ell,\ell
'})_{s\in\spe,\ell,\ell'\geq1}$, we proved that the quantities
%
\begin{equation}
\frac{1}{k} \biggl( \e\log\sum_{\alpha\in\Natural^r} \prod
_{i\in I^+} 2\smsp\ch \bigl(C_{i}(\alpha)
\bigr) v_\alpha - \e\log\sum_{\alpha\in\Natural^r} \exp
\bigl(\sqrt{k} D(\alpha ) \bigr) v_\alpha \biggr) \label{lowerend}
\end{equation}
can be used to approximate the lower limit of the free energy. It
remains to observe that, similarly to (\ref{simp1}) and (\ref{simp2}),
using standard properties of the Ruelle probability cascades (again, we
refer to the proof of Lemma~3.1 in \cite{SKmodel}),
\begin{eqnarray*}
\frac{1}{k} \e\log\sum_{\alpha\in\Natural^r} \prod
_{i\in I^+} 2\smsp\ch \bigl(C_{i}(\alpha) \bigr)
v_\alpha &= &\frac{1}{k} \sum_{i\in I^+}
\e\log\sum_{\alpha\in\Natural^r} 2\smsp\ch \bigl(C_{i}(
\alpha) \bigr) v_\alpha
\\
&=& \sum_{s\in\spe} \lambda_s \e\log\sum
_{\alpha\in\Natural^r} 2\smsp\ch \bigl(C^s(\alpha) \bigr)
v_\alpha
\\
&=& \log2+ \sum_{s\in\spe} \lambda_s
X^s_0
\end{eqnarray*}
and
\[
\frac{1}{k} \e\log\sum_{\alpha\in\Natural^r} \exp \bigl(
\sqrt{k} D(\alpha ) \bigr) v_\alpha = \frac{1}{2} \sum
_{0\leq\ell\leq r-1} \zeta_\ell (Q_{\ell+1} -
Q_{\ell} ).
\]
Therefore, (\ref{lowerend}) is precisely $\PP(\zeta,q)$ defined in
(\ref
{Pzeta}), and this finishes the proof of the lower bound.

\section{Some open questions}

An obvious question that arises is what happens when $\Delta^2$ is not
positive definite, for example, in the case of a bipartite model with
two interacting species and no interactions within species, that is,
$\Delta_{12}^2>0$, and $\Delta_{11}^2 = \Delta_{22}^2=0$. Notice that
our proof of the lower bound for the free energy still works in this
case, but the Guerra-type upper bound in Section~\ref{Sec2label}
utilized the condition $\Delta^2\geq0$ in an essential way.

It was clear from the proof of the lower bound, in particular from the
equations~(\ref{zetafop}) and (\ref{zetafops}), that the parameters
$(\zeta_\ell)$ in (\ref{zetas}) and $(q_\ell)$ in (\ref{qs}) can be
interpreted as encoding the joint distribution of the overlaps within
species and, therefore, the minimizer in formula (\ref{Parisi}) for the
free energy has an important physical interpretation. As a result, as
in the original Sherrington--Kirkpatrick model, there are many
interesting questions about this formula that one can study. For
example, can one extend the result in \cite{AC} to show the uniqueness
of this minimizer? The main result in \cite{AC} implies that the
functional in (\ref{Pzeta}) is strictly convex in the vector $(\zeta
_\ell)_{\ell\leq r}$ for fixed parameters (\ref{qs}), which is
sufficient to prove the uniqueness of the minimizer for one system, but
not obviously for the multi-species case. Another important problem
would be to understand the phase transition in this model and to
describe the replica symmetric (RS) region when the minimizer
corresponds to a distribution (\ref{zetafop}) concentrated on one point
$q\in[0,1]$, that is,
\[
\zeta_0 = 0,\qquad \zeta_1=1,\qquad \zeta_2=1,\qquad
q_0=0, \qquad q_1 = q,\qquad q_2=1.
\]
For technical reasons (to define the Ruelle probability cascades) we
assumed that the inequalities in (\ref{zetas}) are strict, but the
infimum in (\ref{Parisi}) may be achieved on the limiting case when
some inequalities become equalities. If the infimum is replica
symmetric, it is easy to write down the following critical point
equations for the parameters $q^s = q_1^s$ for $s\in\spe$:
\[
\sum_{s\in\spe} \lambda_s
\Delta_{st}^2 \bigl( q^s - \e
\mbox{th}^2 \bigl(z \sqrt{Q^s} + h_s \bigr)
\bigr) = 0 \qquad\mbox{for all } t\in\spe,
\]
where $z$ is a standard Gaussian random variable, $Q^s = 2\sum_{t\in
\spe} \Delta_{st}^2 \lambda_t q^t$ and $(h_s)_{s\in\spe}$ is a vector
of external fields corresponding to each species. (For simplicity of
notation, we did not consider external fields above, but including them
does not affect any arguments.) Assuming that $\Delta^2$ is invertible,
this system is equivalent to
\[
q^s = \e \mbox{th}^2 \bigl(z \sqrt{Q^s} +
h_s \bigr) \qquad\mbox{for all } t\in\spe.
\]
In the SK model, this reduces to one equation, and the uniqueness of
its solution is known as the Latala--Guerra lemma; see Section A.14 in
\cite{SG2}. It would be interesting to see if the solution of the above
system of equations is also unique. In that case, it should not be
difficult to prove replica symmetry breaking above some analogue of the
AT line (in this case, some surface) by the same method as in the SK
model; see \cite{TAT} or Theorem~13.3.1 in \cite{SG2}. However, to
characterize the replica symmetric region exactly, one would probably
need to work much harder. Notice that the multi-species model allows
for some interesting possibilities; for example, one can imagine that
for some choice of parameters, the replica symmetry is broken in some
species but not the others.

\section*{Acknowledgments} The author would like to thank Wei-Kuo Chen
for several helpful discussions and the referee for several comments
about the paper.

%

%




\printaddresses
\end{document}